\documentclass[a4paper,twoside,11pt]{article}
\usepackage{amssymb}
\usepackage{amsmath}
\usepackage{amsthm}
\usepackage{latexsym}
\usepackage{amsfonts}
\usepackage{graphicx}
\usepackage{graphics}

\UseRawInputEncoding
\newcommand{\dis}{\displaystyle}
\theoremstyle{plain}
\newtheorem{thm}{Theorem}[section] 
\newtheorem*{thm*}{Theorem}
\newtheorem{prop}{Proposition}[thm]

\newtheorem*{clm*}{Claim}

\newtheorem{Def}{Definition}[thm]
\newtheorem*{term}{Terminology}
\newtheorem*{Proofclm}{Proof of Claim}

\theoremstyle{definition}
\newtheorem{rem}[thm]{Remark}

\newtheorem*{Proof}{Proof}

\newcommand{\el}{\ell}

\newcommand{\ra}{\;\rightarrow\;}

\newcommand{\de}{\delta }

\newcommand{\thi}{\theta }

\newcommand{\La} {{\mit\Lambda}}

\newcommand{\la}{\lambda }

\newcommand{\si}{\sigma }

\newcommand{\C}{\mathbb{C}}

\newcommand{\R}{\mathbb{R}}

\newcommand{\ssum}{\sum\limits}

\newcommand{\ld}{\ldots}

\newcommand{\qb}{$\quad\blacksquare$}

\begin{document}
\pagestyle{myheadings}
\markboth{Some results of topological genericity}{C. Pandis}
\title{\bf Some results of topological genericity }
\author{C. Pandis}
\date{}
\maketitle
\begin{abstract}
We show topological genericity for the set of functions in the space $\bigcap\limits_{p<1}H^p$ on the open unit disc such that the sequence of Taylor coefficients of the function and of all derivatives of the function  are unbounded. Results of similar nature are valid when the space $\bigcap\limits_{p<1}H^p$ is replaced by 
$H^p(0<p<1)$ and by localized versions of such spaces. Looking at the smaller  space $A(\mathbb{D}) \subseteq H^{\infty}$ we show topological genericity for the set of functions in $A(\mathbb{D})$ and of all derivatives such that the sequence of Taylor coefficients of the function are outside of  $\el^1$.

We also show topological genericity for the set of functions in the space $\bigcap\limits_{p<1}h^p$ whose  harmonic conjugate  does not belong in any $h^q(q>0)$
\end{abstract}
\smallskip\\
{\em AMS classification numbers}: 30H10

\section{Introduction}\label{sec1}
\noindent
  If there exists an object with a ``bad'' property, then a general principle is that there are many such objects and their set is big in various senses. In this paper we study such properties. It is well known that for a function  $f$ in $H^1$ its Taylor coefficients tend to zero due to Rieman-Lebesgue Theorem. The space  $\bigcap\limits_{p<1}H^p$, which is very close to $H^1$  has a diametrically opposite property; a generic function of this space has Taylor coefficients that are unbounded. More generally, a generic function of $\bigcap\limits_{p<1}H^p$ has derivatives with unbounded Taylor coefficients. Those results can also be generalised for the local spaces $H^p_{[A,B]}$ and $\bigcap\limits_{p<1}H^p_{[A,B]}$. Taking a step further we show that a generic function and its derivatives in the  disc algebra $A(\mathbb{D})$ has Taylor coefficients outside of $\el^1$.

In the second part  we prove a generic result concerning the harmonic hardy spaces $h^p$. M.Riez proved that if $1<p<\infty$ the space $h^p$ is ''self-conjugate'' in the sense that if $u\in h^p$ then the same is true for the harmonic conjugate $\tilde{u}$. This is not true for $h^1$, but here Kolmogorov's theorem provides a substitute : if $u$ in $h^1$ then $\widetilde{u}$ is in $ h^p $ for all $p<1$ . We show that in the space $\bigcap\limits_{p<1}h^p$, a generic function $u\in \bigcap\limits_{p<1}h^p$ satisfies 
\[ 
 \sup_{0<r<1}\dfrac{1}{ 2\pi}\int_{-\pi}^{\pi}|\widetilde{u}(re^{i\thi})|^q d\thi = +\infty \quad \text{for all} \quad q>0 
\]  

\section{Preliminaries}\label{sec:2}
\noindent

Let $\mathbb{D} =\{z\in\C:|z|<1\}$ be the open unit disc. A holomorphic function $f:\mathbb{D}\rightarrow \C$ belongs to the Hardy space $H^p$ $(0<p<+\infty)$ if $\dis\sup_{0<r<1}\dfrac{1}{2\pi}\int^{2\pi}_0|f(re^{i\thi})|^pd\thi<+\infty$. It belongs to the Hardy space $H^\infty$ if $\dis\sup_{|z|<1}|f(z)|<+\infty$. The space $H^\infty$ endowed with the supremum norm on $\mathbb{D}$ is a Banach space, but polynomials are not dense in this space. For $1\leq p<+\infty$ the space $H^p$ endowed with the norm
\[\|f\|_p=\sup_{0<r<1}\bigg\{\frac{1}{2\pi}\int^{2\pi}_0|f(re^{i\thi})|^pd\thi\bigg\}^{1/p}
\]
is also a Banach space. Then for $f,g\in H^p$ their distance is
\[
d_p(f,g)=\sup_{0<r<1}\bigg\{\frac{1}{2\pi}\int^{2\pi}_0|f(re^{i\thi})-g(re^{i\thi})|^p
d\thi\bigg\}^{1/p}, \ \ 1\le p<+\infty.
\]
For $0<p<1$ we endow $H^p$ with the metric
\[
d_p(f,g)=\sup_{0<r<1}\frac{1}{2\pi}\int^{2\pi}_0|f(re^{i\thi})-g(re^{i\thi})|^pd\thi, \ \ f,g\in H^p, \ \ 0<p<1
\]
and then $H^p$ becomes a topological vector space endowed with a metric invariant by translation which is complete ($F$-space). For $0<p<+\infty$ polynomials are dense in $H^p$. Also convergence in $H^p$, $0<p\le+\infty$ implies uniform convergence on each compact subset of $\mathbb{D}$.
For $a,b\in(0,+\infty]$, $a<b$ we have $H^b\subset H^a$ and the injection map is \linebreak continuous. Jensen's inequality implies that the map
\[
a\ra\sup_{0<r<1}\Big\{\frac{1}{2\pi}\int^{2\pi}_0|f(re^{i\thi})|^ad\thi\Big\}^{1/a}
\]
is increasing. Obviously we also have
\[
\sup_{0<r<1}\Big\{\frac{1}{2\pi}\int^{2\pi}_0|f(re^{i\thi})|^pd\thi\Big\}^{1/p}\le
\sup_{|z|<1}|f(z)|.
\]

Next for $0<a\le+\infty$ we consider the intersection $\bigcap\limits_{p<a}H^p$.

Convergence in this space is equivalent with convergence in all spaces $H^p$, $p<a$. Equivalently, we consider a strictly increasing sequence $p_n$ converging to $a$ and the metric in $\bigcap\limits_{p<a}H^p$ is defined by
\[
d(f,g)=\sum^\infty_{n=1}\frac{1}{2^n}\frac{d_{p_n}(f,g)}{1+d_{p_n}(f,g)}, \ \ f,g\in\bigcap_{p<a}H^p.
\]
This space is also complete, in fact an $F$-space. Obviously convergence in $\bigcap\limits_{p<a}H^p$ implies uniform convergence on each compact subset of $\mathbb{D}$.
\begin{prop}\label{prop2.1}
Polynomials are dense in $\bigcap\limits_{p<a}H^p$, for every $a$, $0<a\le+\infty$.
\end{prop}

For the proof it suffices for $f\in\bigcap\limits_{p<a}H^p$ to control $d_{p_n}(f,P)$, for $n=1,\ld,N$ for any finite $N$. Because of the monotonicity of the map $p\ra\dis\sup_{0<r<1}\Big\{\dfrac{1}{2\pi}\int^{2\pi}_0|f(re^{i\thi})-P(re^{i\thi})|^pd\thi\Big\}^{1/p}$
is suffices to control $d_{p_N}(f,P)$. But this is possible, because polynomials are dense in $H^{p_N}$ since $p_N<+\infty$. Thus Proposition
\ref{prop2.1} holds.

Next we present localized versions of the previous spaces.

Let $0<p<+\infty$ and $A,B\in\R$, $A<B$. Then a holomorphic function $\linebreak$  $f:\mathbb{D} \rightarrow \C$ belongs to $H^p_{[A,B]}$ if $\dis\sup_{0<r<1}\int^B_A|f(re^{i\thi})|^p\dfrac{d\thi}{B-A}<+\infty$ and to $H^\infty_{[A,B]}$ if $\dis\sup_{0<r<1}\dis\sup_{A\le\thi\le B}|f(re^{i\thi})|<+\infty$. Because of the monotonicity of the function \[ a\ra\dis\sup_{0<r<1}\Big\{\int^B_A|f(re^{i\thi})|^a\dfrac{d\thi}{B-A}\Big\}^{1/a}\] it follows $H^b_{[A,B]}\subset H^a_{[A,B]}$ for  $0<a<b\le+\infty$.

Convergence in $H^p_{[A,B]}$ of a sequence $f_n$ towards $f$ where $f_n,f\in H^p_{[A,B]}$ is equivalent to uniform convergence on all compact subsets of $\mathbb{D}$ and
\[
\sup_{0<r<1}\bigg\{\int^B_A|f_n(re^{i\thi})-f(re^{i\thi})|^p\frac{d\thi}{B-A}\bigg\}^{1/p}
\xrightarrow{n\ra +\infty}{}0 \ \ \mbox{for} \ \ 0<p<+\infty \ \ \mbox{and}
\]
\[
\sup_{0<r<1}\sup_{A\le\thi\le B}|f_n(re^{i\thi})-f(re^{i\thi})|\xrightarrow{n\ra +\infty}{}0 \ \ \mbox{for} \ \ p=+\infty.
\]
The metric giving this topology in $H^p_{[A,B]}$ is defined by
\begin{align*}
d_{p,[A,B]}(f,g)=&\sup_{0<r<1}\bigg\{\int^B_A|f_n(re^{i\thi})-g(re^{i\thi})|^p\cdot
\frac{d\thi}{B-A}\bigg\}^{1/p}\\
&+\sum^\infty_{n=2}\frac{1}{2^n}\frac{\dis\sup_{|z|\le1-\frac{1}{n}}|f(z)-g(z)|}
{1+\dis\sup_{|z|\le1-\frac{1}{n}}|f(z)-g(z)|} \ \ \mbox{for} \ \ 1\le p<+\infty
\end{align*}
\begin{align*}
d_{p,[A,B]}(f,g)=&\sup_{0<r<1}\int^B_A|f(re^{i\thi})-g(re^{i\thi})|^p\frac{d\thi}{B-A} \\
&+\sum^\infty_{n=2}\frac{1}{2^n}\frac{\dis\sup_{|z|\le1-\frac{1}{n}}|f(z)-g(z)|}
{1+\dis\sup_{|z|\le1-\frac{1}{n}}|f(z)-g(z)|} \ \ \mbox{for} \ \ 0<p<1 \ \ \mbox{and}
\end{align*}
\begin{align*}
d_{\infty,[A,B]}(f,g)=&\sup_{0<r<1}\sup_{A\le\thi\le B}|f(z)-g(z)| \\
&+\sum^\infty_{n=2}\frac{1}{2^n}\frac{\dis\sup_{|z|\le1-\frac{1}{n}}|f(z)-g(z)|}
{1+\dis\sup_{|z|\le1-\frac{1}{n}}|f(z)-g(z)|} \ \ \mbox{for} \ \ p=+\infty.
\end{align*}
Obviously, convergence in $H^p_{[A,B]}$ implies uniform convergence on all compact subsets of $\mathbb{D}$. Also for $0<a<b\le+\infty$ the injection map $H^b_{[A,B]}\subset H^a_{[A,B]}$ is continuous. Finally these spaces are complete, in fact $F$-spaces and $H^p_{[A,B]}=H^p$ when $B-A\ge 2\pi$ and in general $H^p\subset H^p_{[A,B]}$, provided $A<B$.

Let $0<a\le+\infty$. Convergence in the space $\bigcap\limits_{p<a}H^p_{[A,B]}$ is equivalent to convergence in all $H^p_{[A,B]}$ for $p<a$. A metric in $\bigcap\limits_{p<a}H^p_{[A,B]}$ compatible with this topology is given by
\[
d(f,g)=\sum^\infty_{n=1}\frac{1}{2^n}\frac{d_{p_n,[A,B]}(f,g)}{1+d_{p_n,[A,B]}(f,g)}
\]
where $p_n$ is any strictly increasing sequence converging to $a$. This space is complete, in fact an $F$-space. Obviously convergence in $\bigcap\limits_{p<a}H^p_{[A,B]}$ implies uniform convergence on all compact subsets of $\mathbb{D}$

A function $f: \overline{\mathbb{D}} \rightarrow \mathbb{C}$ belongs to the  disc algebra  $A(\mathbb{D})$ if it is holomorphic on the open unit disc and continuous on the closed unit disc or else we can say $A(\mathbb{D})= H^{\infty}(\mathbb{D})\cap C(\overline{\mathbb{D}})$. Endowed with the unifrom norm $A(\mathbb{D})$ becomes a  Banach space (a Banach algebra in particular). Obviously convergence in $A(\mathbb{D})$ implies uniform convergence on $\overline{\mathbb{D}}$ which in turn implies uniform convergence on compact subsets of $\mathbb{D}$. Polynomials are dense in $A(\mathbb{D})$.

A harmonic function $u:\mathbb{D}\!\ra\!\R$  belongs in  $h^p$  if   $\dis\sup_{0<r<1}\dfrac{1}{2\pi}\int^{\pi}_{-\pi}|u(re^{i\thi})|^pd\thi<+\infty$.
For $p\geq 1$ then space $h^p$ endowed with the norm 
\[ \|u\|_p =\dis\sup_{0<r<1}\bigg\{\dfrac{1}{2\pi}\int^{\pi}_{-\pi}|u(re^{i\thi})|^pd\thi \bigg\}^{1/p}  \]
is a Banach space.

For $0<p<1$ we endow $h^p$ with the translation invariant metric
\[  
d_p(u,v)=  \sup_{0<r<1}\dfrac{1}{2\pi}\int^{\pi}_{-\pi}|u(re^{i\thi})-v(re^{i\thi})|^pd\thi \]
 In either case metric convergence implies uniform convergence in compact subsets of $\mathbb{D}$, so even if $0<p<1$, the space $h^p$ is complete. For $a,b\in(0,+\infty]$, $a<b$ we have $h^b\subset h^a$ as before.\cite{8}

 Harmonic cojugate is only defined to within an additive constant, working at the unit cirlce, it is customary to require that its value at zero is zero. $\cite{2}$

Next for $0<a < 1 \ $ we consider the intersection $\bigcap\limits_{p<a}h^p$.

Convergence in this space is equivalent with convergence in all spaces $h^p$, $p<a$. Equivalently, we consider a strictly increasing sequence $p_n$ converging to $a$ and the metric in $\bigcap\limits_{p<a}h^p$ is defined by
\[
d(u,v)=\sum^\infty_{n=1}\frac{1}{2^n}\frac{d_{p_n}(u,v)}{1+d_{p_n}(u,v)}, \ \ u,v\in\bigcap_{p<a}h^p.
\]
This space is also complete . Obviously convergence in $\bigcap\limits_{p<a}h^p$ implies uniform convergence on each compact subset of $\mathbb{D}$.

We now mention a theorem of Borel and Carath\'eodory \cite{9} that shows that an analytic function may be bounded by its real part that will be of use later in the paper.

\begin{thm*}
    Let a function $f$  be analytic on a closed disc of radius R, centered at the origin. Suppose that $ r < R$. Then, we have the following inequality:
\[ 
\sup_{|z|\leq r }|f(z)| \leq \dfrac{2r}{R-r} \sup_{|z|\leq R}\mathfrak{R}(f(z)) + \dfrac{R+r}{R-r}|f(0)|
\]

\end{thm*}

In order to show that the sets of these functions are ``big'' we  prove a slightly changed variation of a result of M. Siskaki (\cite{4}), that was introduced later by  V.Nestoridis and E.Thirios (\cite{6}).

\begin{prop}[V.Nestoridis, E.Thirios variation]
Let $V$ be a topological vector space over the field $\R$ or $\C$. Let $X$ be a non empty set and $\C^X$ the set of complex functions defined on $X$. Let $T:V\ra\C^X$ be such that

1) For every $x\in X$ the function $V\ni f\ra T(f)(x)\in\C$ is continuous.

2) $|T(f-g)(x)|\le|T(f)(x)|+|T(g)(x)|$ for all $f,g\in V$ and $x\in X$.

3) For every $f\in V$, if $T(f)$ is unbounded on $X$, then there is a sequence $(\la_n)_n$ of numbers in $\R$ or in $\C$, respectively, with $\la_n\ra0$ as $n\ra\infty$ such that $T(\la_nf)$ is unbounded on $X$ for every $n\ge1$.

We set $S=\{f\in V:T(f)$ is unbounded on $X\}$. Then, either $S=\emptyset$ or $S$ is a $G_\de$ and dense subset of $V$.
\end{prop}

For the needs of this paper we will notice that condition $2$ of this variation can be replaced as below 

\begin{prop}\label{prop2.5}
Let $V$ be a topological vector space over the field $\R$ or $\C$. Let $X$ be a non empty set and $\C^X$ the set of complex functions defined on $X$. Let $T:V\ra\C^X$ be such that

1) For every $x\in X$ the function $V\ni f\ra T(f)(x)\in\C$ is continuous.

2) $|T(f-g)(x)|\le M(|T(f)(x)|+|T(g)(x)|)$ for all $f,g\in V$, for some $M>0$ and $x\in X$.

3) For every $f\in V$, if $T(f)$ is unbounded on $X$, then there is a sequence $(\la_n)_n$ of numbers in $\R$ or in $\C$, respectively, with $\la_n\ra0$ as $n\ra\infty$ such that $T(\la_nf)$ is unbounded on $X$ for every $n\ge1$.

We set $S=\{f\in V:T(f)$ is unbounded on $X\}$. Then, either $S=\emptyset$ or $S$ is a $G_\de$ and dense subset of $V$.
\end{prop}

\begin{Proof}
The proof that $S$ is a $G_\de$ is omitted, because it follows simply from 1 and is similar to the proof in \cite{4}.

Suppose $S\neq\emptyset$ and let $f\in S$. Thus, $T(f)$ is unbounded on $X$. If $S$ is not dense, then there exist $g\in V$ so that $g\notin\overline{S}$. Then $T(g)$ is bounded on $X$ by a constant $M_1$.

Since $V$ is a topological vector space it holds $g+\la_nf\ra g$ as $n\ra+\infty$.

According to our assumptions we have
\begin{align*}
|T(\la_nf)(x)|=|T(g+\la_nf-g)(x)|&\le M|T(g+\la_nf)(x)|+M|T(g)(x)| \\
&\le M|T(g+\la_nf)(x)|+MM_1
\end{align*}
for all $x\in X$ where $M_1<+\infty$ and $M>0$ is independent of $x$.\\
Then for every $n\ge1$ it follows that $T(g+\la_nf)$ is unbounded on $X$ and it is deduced that $g$ belongs to the closure of $S$, which is a contradiction.
\qb
\end{Proof}
\section{The Results}\label{sec:3}
\begin{term}
    
From now on we refer as $\beta(f)$ to be the sequence of Taylor \\
coefficients of a holomorphic function f on the unit disc  :  
\[ 
\beta(f) = ( \beta_k (f) )_{k=0}^{\infty}
\]
\end{term}

\begin{thm} 
The set  $\mathcal{A} = \big  \{ f \in \bigcap\limits_{p<1}H^p : \beta(f) \notin \el^{\infty} \big  \}$  is  a $G_{\delta}$ and  dense   subset of $\bigcap\limits_{p<1}H^p$ .
 
\end{thm}
\begin{Proof} 
Let $f(z)= \dfrac{1}{1-z}\log\bigg(\dfrac{1}{1-z}\bigg)$ , $z \in \mathbb{D} $ .Then $f$ belongs to $H^p $ for all $0<p<1$ and its Taylor coefficients are unbounded \cite{10} and thus $\mathcal{A} \neq \emptyset $.

We now notice that if  $P$ is a polynomial and $f \in  \mathcal{A} $  then  $f+P \in \mathcal{A}$. Since the set of polynomials is dense in $\bigcap\limits_{p<1}H^p $, it follows that the set $\{f+P:P$\, polynomial$\}$ is dense in  $\bigcap\limits_{p<1}H^p$. Since the last set is contained in $\mathcal{A}$, it follows that $\mathcal{A}$ is dense in $\bigcap\limits_{p<1}H^p$.

In order to show that $\mathcal{A}$ is $G_{\delta}$ it suffices to prove that $\bigcap\limits_{p<1}H^p \setminus \mathcal{A}$   is a denumerable union of closed subsets of $ \bigcap\limits_{p<1}H^p$.

For $M$ a natural number we consider the set  
\[ 
\Omega_M= \bigg \{f \in \bigcap\limits_{p<1}H^p : |\beta_n(f)| \leq M  \;\; , \forall n \in \mathbb{N} \bigg \}
\]
Then $ \bigcap\limits_{p<1}H^p \setminus \mathcal{A}= \bigcup\limits_{M} \Omega_M$ . 

We verify that each set $\Omega_M$ is a closed subset of $\bigcap\limits_{p<1}H^p$. Indeed, let $f_k \in \Omega_ M$ be a sequence converging in $\bigcap\limits_{p<1}H^p$ to some $f \in \bigcap\limits_{p<1}H^p$. Then $f_k$ converges uniformly to f  on  compacta   of $\mathbb{D}$ ,which implies $\beta_n(f_k) \xrightarrow[k\to\infty]{}\beta_n(f) , \forall n \in \mathbb{N}$ .

Since $|\beta_n(f)|=\lim_{k \to \infty}|\beta_n(f_k)| \leq M$ and thus $f\in \Omega_M$ .We shown that $\bigcap\limits_{p<1}H^p \setminus \mathcal{A}$ is $F_{\sigma}$  which completes the proof .  $\blacksquare$

\end{Proof}

\begin{thm}
The sets  $\mathcal{C}_n =  \big \{ f \in \bigcap\limits_{p<1}H^p : \beta(f^{(n)}) \notin \el^{\infty} \big \} $ 
, $(n \in \mathbb{N})$ are  $G_{\delta}$ and  dense   subsets of $\bigcap\limits_{p<1}H^p$.
\end{thm}
\begin{Proof}
 It is obvious that  $\beta_k(f^{(n)})=(k+n)(k+(n-1))\cdots(k+1)\beta_{k+n}(f)$.

We  notice that for a function $f \in \bigcap\limits_{p<1}H^p $ such that the sequence $\beta(f)$ is unbounded then the same is true for $\beta(f^{(n)})$ and so $ \mathcal{A} \subseteq \mathcal{C}_n$ for all $n \in \mathbb{N}$. Since $\mathcal{A}$ is non-empty and dense subset of $\bigcap\limits_{p<1}H^p$, so is $\mathcal{C}_n$  for all $n \in \mathbb{N}$.

We fix a natural number $n$ and proceed to show that  $\bigcap\limits_{p<1}H^p \setminus \mathcal{C}_n$ is a denumerable union of  closed 
 subsets of $\bigcap\limits_{p<1}H^p $.

 As before, for $M$ a natural number we consider the set 
 \[
    \Omega_{M}^n= \bigg \{f \in\bigcap\limits_{p<1}H^p :  |\beta_k(f^{(n)})| \leq M  \;\; , \forall k \in \mathbb{N} \bigg \}
\]
and notice that  $\bigcap\limits_{p<1}H^p \setminus \mathcal{C}_n=\bigcup\limits_{M} \Omega_{M}^n$. 

The set $\Omega_ {M}^n$ is closed, indeed let $f_l \in  \Omega_ {M}^n$ be a sequence  converging in $\bigcap\limits_{p<1}H^p$ to some $f \in \bigcap\limits_{p<1}H^p$. Then $f_l \rightarrow f$ uniformly on  compacta   of $\mathbb{D}$ which implies, from Weierstrass theorem, that $f_l^{(n)}$ converges uniformly to $f^{(n)}$ on compacta of $\mathbb{D}$, for all $n \in \mathbb{N}$, and thus  $\beta_k(f_l^{(n)}) \xrightarrow[l\to\infty]{}\beta_k(f^{(n)})$ ,  for all  $ n \in \mathbb{N}$ .

Since $|\beta_k(f^{(n)})|=\lim_{l \to \infty}|\beta_k(f_{l}^{(n)})| \leq M$, for all $k\in \mathbb{N}$ we obtain  $f \in \Omega_{M}^n$. We shown that $\bigcap\limits_{p<1}H^p \setminus \mathcal{C}_n$ is $F_{\sigma}$   $\blacksquare$

\end{Proof}

\begin{rem}
The above results holds true if we replace the space $\bigcap\limits_{p<1}H^p$ with each of the spaces  $H^p$ for $0<p<1$, or even with their local spaces $H^p_{[A,B]}$ for $0<p<1$ and $\bigcap\limits_{p<1}H^p_{[A,B]}$, since $\bigcap\limits_{p<1}H^p \subseteq H^p \subseteq H^p_{[A,B]} $ and so $\dfrac{1}{1-z}\log\bigg(\dfrac{1}{1-z}\bigg) \in H^p_{[A,B]}$ for all $0<p<1$.

\end{rem}

 \begin{rem}
In contrast with the situation of the Taylor coefficients of the generic function $f$ and its derivative, the same is not true for the primative $F(f)$ of the function. We see that if $f \in  \bigcap\limits_{p<1}H^p $ then $\beta(F(f))$   is bounded  and that is immediate by using a variation of a known theorem  asserting that if  $f\in \bigcap\limits_{p<1}H^p $ then $F(f) \in \bigcap\limits_{p< +\infty}H^p \subseteq H^1$ \cite{1}, \cite{6}.

 Here it would be interesting to mention a result of V.Nestoridis \cite{7}, that for a generic function $f$ on $H^1$ the sequence $ \beta(F(f)) $ is outside any $\el^p$ space smaller than $\el^1$ i.e. with $0<p<1$; thus, $\beta(F(f)) \in\el^1\setminus \Big(\dis\bigcup_{0<p<1}\el^p\Big)$ holds generically for every $f$ in $H^1$. $\blacksquare$
 \end{rem} 

We now notice that for a generic function $f\in A(\mathbb{D}) \subseteq H^{\infty}  $, then the weaker condition $\beta(f)\notin \el ^1$ holds.

A direct consequence of Hardy's inequality \cite{1} is that if $f'\in H^1$ then $\beta(f) \in \el^1$. In particular
$\beta(f) \in \el^1$ if $f$ is a conformal mapping of the unit disc onto a Jordan domain with rectifiable boundary.\cite{1}

Golubev in \cite{11} gave the first example of a conformal mapping of the unit disc $\mathbb{D}$ onto a Jordan domain which carries a boundary set of measure zero onto a set of positive measure on the nonrectifiable boundary of the image domain .
\begin{thm}
    The set  $\mathcal{E} = \big  \{ f \in A(\mathbb{D}) : \beta(f) \notin \el^{1} \big  \}$  is  a $G_{\delta}$ and  dense   subset of $A(\mathbb{D})$.
\end{thm}

\begin{Proof}
    We mentioned before that $\mathcal{E} \neq\emptyset$, let $g \in \mathcal{E}$. To prove that $\mathcal{E}$ is a dense subset of $A(\mathbb{D})$ we will show that the interior of $A(\mathbb{D})\setminus \mathcal{E}$ is void. Assume $f\in (A(\mathbb{D})\setminus \mathcal{E})^{\mathrm{o}}$ then 
    \[
    f +\frac{1}{n}g \xrightarrow[n \ra \infty]{\|\cdot \|_{\infty}} f 
    \]
and $f$ is in the interior of $A(\mathbb{D})\setminus \mathcal{E}$. It follows that that for some $n_0\in \{1,2,\ldots \}$ the function $f+\frac{1}{n_0}g$ belongs to $A(\mathbb{D})\setminus \mathcal{E}$.
This suggests that $\beta(\frac{1}{n_0}g) \in \el^1$ which is a contradiction. Thus $(A(\mathbb{D})\setminus \mathcal{E})^{\mathrm{o}}=\emptyset$

In order to show that $\mathcal{E}$ is a $G_\de$ is suffices to prove that $A(\mathbb{D})\setminus\mathcal{E}$ is a denumerable union of closed subsets of $A(\mathbb{D})$.

For $M$ and $N$ natural numbers we consider the set
\[
\Omega_{M,N}=\bigg\{f\in A(\mathbb{D}) :f(z)=\sum^\infty_{n=0}\beta_n(f)z^n, \;\;\sum^N_{n=0}|\beta_n(f)|\le M\bigg\}.
\]
Then $A(\mathbb{D}) \setminus \mathcal{E}= \dis\bigcup_M\Big[\dis\bigcap_N \Omega_{M,N }\Big]$.

We verify that each set $\Omega_{M,N}$ is a closed subset of $A(\mathbb{D})$. Indeed, let $f_k\in \Omega_{M,N}$ be a sequence converging in $A(\mathbb{D})$ to some $f\in A(\mathbb{D})$. Then $f_k$ converges uniformly on compacta of $\mathbb{D}$ to $f$, which implies $\beta_n(f_k)\xrightarrow[k\ra\infty]{}\beta_n(f)$ for every $n=0,1,2,\ld$ . Since $\ssum^N_{n=0}|\beta_n(f_k)|\le M$ for all $k$, it follows $\ssum^N_{n=0}|\beta_n(f)|\le M$; that is, $f\in \Omega_{M,N}$ and the set $\Omega_{M,N}$ is closed in $A(\mathbb{D})$. The same holds for the intersections $\dis\bigcap_N \Omega_{M,N}$ and their denumerable union $A(\mathbb{D})\setminus\mathcal{E} $ is an $F_\si$.
The proof is complete. $\blacksquare$
\end{Proof}
 
\begin{thm}
    The sets $\mathcal{R}_n =\big  \{ f \in A(\mathbb{D}) : \beta(f^{(n)}) \notin \el^{1} \big  \}$ $(n>1)$ are   $G_{\delta}$ and  dense   subsets of $A(\mathbb{D})$.
\end{thm}

\begin{Proof}
    Since  $\beta_k(f^{(n)})=(k+n)(k+(n-1))\cdots(k+1)\beta_{k+n}(f)$, it is immidiate that $\mathcal{E}\subset \mathcal{R}_n$ for  $n>1$ and $\mathcal{E}$ is a dense subset of $A(\mathbb{D})$ then so is $\mathcal{R}_n$.

As before
for $M$ and $N$ natural numbers we consider the set
\[
\Omega_{M,N}^n=\bigg\{f\in A(\mathbb{D}) :f(z)=\sum^\infty_{k=0}\beta_k(f)z^k, \;\;\sum^N_{k=0}|\beta_k(f^{(n)})|\le M\bigg\}.
\]
Then $A(\mathbb{D}) \setminus \mathcal{R}_n= \dis\bigcup_M\Big[\dis\bigcap_N \Omega_{M,N }^n \Big]$.

Now we show that the set  $A(\mathbb{D}) \setminus \mathcal{R}_n$ is $F_{\sigma}$.

Let $f_l\in \Omega_{M,N}^n $ be a sequence converging in $A(\mathbb{D})$ to some $f\in A(\mathbb{D})$. Then $f_l$ converges unifromly on compacta of $\mathbb{D}$ to $f$ again by Weierstrass theorem we obtain  $\beta_k(f_{l}^{(n)})\xrightarrow[l\ra\infty]{}\beta_k(f^{(n)})$. As seen before the set $\Omega_{M,N}^n$  is closed. $\blacksquare$
\end{Proof}

We now use Baire's category theorem and the variation of Siskaki's proposition  to prove that for a generic function $ u $ in  $\bigcap\limits_{p<1}h^p $ its harmonic conjugate $\widetilde{u}$ is outside of any harmonic hardy space $h^q$  for   $q>0$ .

We may mention here again that the harmonic conjugate is normalised i.e its value at zero is zero.
\begin{Def}
 We denote by $\La_q$, ($q>0$) the set of all functions $u\in \bigcap\limits_{p<1}h^p$ such that the  harmonic conjugate $\widetilde{u} \notin h^q $.

\end{Def}

\begin{prop} The sets $\La_q$ for $0<q<1$ are  $G_{\delta} $ and dense subset of $\bigcap\limits_{p<1}h^p$ .

\end{prop}

\begin{Proof}
Let 
\[  
    f(z)=u(z)+iv(z)=\sum_{n=1}^{+\infty}\epsilon_n \dfrac{z^{2^n}}{1-z^{2^{n+1}}} ,\quad  \epsilon_n \in \{+1,-1 \}
\]
Then for every choice of signs $\epsilon_n$, $u \in h^p $ for all $p<1$, while for almost every  sequence of signs, $f(z)$ has a radial limit on no set of positive measure. In particular, some choice of the $\epsilon_n$ gives a function $f$ which is not even of Nevanlinna class $N$, but whose real part belongs to $h^p$ for all $p<1$. \cite{1} , and thus $\La_q \neq \emptyset  $ for all $q>0$.

We consider the operator $ T: \bigcap\limits_{p<1}h^p \times (0,1) \ra \C $ given by 
 
\[
T(u,r)=  \frac{1}{ 2\pi}\int_{-\pi}^{\pi}|\widetilde{u}(re^{i\thi})|^q d\thi ,\quad u \in \bigcap\limits_{p<1}h^p  ,  \quad  r\in (0,1) , \, \quad \text{for} \quad  0<q<1
\]

and notice that $\La_q=  \{ u \in \bigcap\limits_{p<1}h^p : \dis\sup_{r\in (0,1)}|T(u,r)|=+\infty \}$.

We are going now to verify the conditions of the variation of M.Siskaki theorem.

\begin{enumerate}
    \item Since $\widetilde{u+v}=\widetilde{u}+\widetilde{v}$ , $\widetilde{\lambda u}=\lambda\widetilde{u}$ and by using the triangular inequality and the well known inequality $|a \pm b|^q \leq 2^q(|a|^q + |b|^q)$ it is immediate that

\[
    |T(u-v,r)| \leq 2^q( |T(u,r)|+|T(v,r)|)
\]
and 
\[ 
|T(\lambda u,r)|= |\lambda|^q|T(u,r)|
\]

\item For all $r\in (0,1)$,  $T( \cdot ,r )$ is continuous :

   \begin{clm*}
Let $(u_m)_{m=1}^{\infty}$ be a sequence of real-valued harmonic functions defined on the unit disc  $\mathbb{D}$ such that $u_m \rightarrow u$   uniformly on compacta of $\mathbb{D}$ then $\widetilde{u_m} \rightarrow \widetilde{u}$  uniformly  on compacta of $\mathbb{D}$.
   \end{clm*}

\begin{Proofclm}
Let $(f_m)_{m=1}^{\infty}$ be a sequence of holomorphic functions on the unit disc $\mathbb{D}$ such that $f_m=u_m+i\widetilde{u_m}$.

We observe that $f_m(0) = u_m(0) \xrightarrow[m\to\infty]{}u(0) = f(0) $.

Let $K$ be a compact subset of the unit disc $\mathbb{D}$, then there exist $0<r<1$ such that $K$ is contained in the open disc $D(0,r)$.

We now pick $R= \dfrac{1+r}{2} $ and so  $D(0,r) \subset D(0,R)\subset \mathbb{D}$, we can now apply the Borel-Carath\'eodory theorem for the holomorphic function \[
f_m-f=u_m-u+i(\widetilde{u_m}-\widetilde{u})
\]
and obtain that
\[  
\sup_{|z|\leq r }|f_m(z)-f(z)| \leq \dfrac{2r}{R-r} \sup_{|z|\leq R}|u_m(z)-u(z)| + \dfrac{R+r}{R-r}|u_m(0)-u(0)|
\] 

Since the disc $ |z| \leq R $ is compact and we know that $u_m\rightarrow u$ uniformly on compact subsets on $\mathbb{D}$ we obtain
\[  
\sup_{z\in K}|f_m(z)-f(z)| \leq \sup_{|z| \leq r }|f_m(z)-f(z)| \xrightarrow[m \to \infty]{} 0
\]
It is now immidiate that $\dis\sup_{z\in K}|\widetilde{u_m}(z) -\widetilde{u}(z)|\xrightarrow[m\to \infty]{}0$. $\blacksquare$

   \end{Proofclm} 
Now let $(u_m)_{m=1}^{\infty}$ be a sequence in $\bigcap\limits_{p<1}h^p$  such that $u_m \rightarrow u$ in  $\bigcap\limits_{p<1}h^p$. Since convergence in $\bigcap\limits_{p<1}h^p$ implies uniform convergence on compacta of $\mathbb{D}$ combined with the above claim, we obtain that $\widetilde{u_m} \rightarrow \widetilde{u}$ uniformly on compacta of $\mathbb{D}$ and so $T(u_m,r) \rightarrow T(u,r)$ for all $r\in(0,1).$
   
\item Is easy to check.

\end{enumerate}

\begin{thm}
The set   $\mathcal{L} =  \big \{ u \in \bigcap\limits_{p<1}h^p : \widetilde{u}\notin h^q, \, \text{for all} \ q>0 \big \} $  is  a $G_{\delta}$ and  dense   subset of $\bigcap\limits_{p<1}h^p$. 
\end{thm}

\begin{Proof}
Because of the monotonicity of the harmonic hardy spaces and that the  assertion that for every $q>0 $ it holds $\widetilde{u} \notin h^q $ is equivalent to the assertion that for every $n\geq 2$ it holds $\widetilde{u} \notin h^{ \frac{1}{n}}$, the set $\mathcal{L}$ is equal to  a denumerable intersection for $n\geq 2 $ of the the sets   $\La_{ \frac{1}{n}}$. According to Proposition 3.1 these sets are $G_{\delta}$ and dense in $\bigcap\limits_{p<1}h^p$. Baire's theorem yields the result. $\blacksquare$

\end{Proof}

\end{Proof}

\begin{rem}
We can replace the space $\bigcap\limits_{p<1}h^p$ with each of the spaces $h^p$ for $0<p<1$.
\end{rem}
\begin{rem}
    One could avoid the use of the variation of M.Siskakis lemma introduced in this paper by noticing that for arbitrary positive numbers $a$ and $b$ and $0<p<1$ it holds
    $(a+b)^p \leq a^p+b^p$ by elementary calculus. \cite{1}
\end{rem}
\begin{rem}
A direct approach on the continuity of the operator without the use of the claim mentioned above  was suggested by V.Nestoridis and A.Siskakis and goes as follows:

Let $0<p<1$ and  $0<r<1$ and let $(u_m)_{m=1}^{\infty}$ be a sequence  such that $u_m \rightarrow u$ in  $\bigcap\limits_{p<1}h^p$. Since convergence in $\bigcap\limits_{p<1}h^p$ implies convergence in each $h^p$ for $p<1$, and so uniform convergence on compacta of $\mathbb{D}$ we get that $u_m \rightarrow u$ uniformly on the disc $r\mathbb{D}$. It is immediate that $u_m \rightarrow u $ in $h^2(r\mathbb{D})$, but in $h^2$  the map sending $u$ to its conjugate $\widetilde{u}$ is an isometry so $\widetilde{u_m} \rightarrow \widetilde{u}$ in $h^2(r\mathbb{D})$, thus $\widetilde{u_m} \rightarrow \widetilde{u} $ in compacta of $r\mathbb{D}$. For every compact $K\subseteq \mathbb{D}$ we can find $0<s<1$ such that $K\subseteq s\mathbb{D}$ and so the result follows.

$\blacksquare$
\end{rem}

\noindent
{\bf Acknowledgements}.
The author would like to express his gratitude towards V.Nestoridis for his interest in this work and his always useful comments and remarks. The author would also like to thank A.Siskakis for the helpful communication.

\vspace*{0.5cm}
\noindent
C. Pandis\\
National and Kapodistrian University of Athens\\
Department of Mathematics\\
Panepistemiopolis, 157 84\\
Athens,\\
Greece\\
e-mail: chrpandis@gmail.com \\
\end{document}